\documentclass[11pt]{amsart}
\usepackage[a4paper,margin=2.6cm]{geometry}
\usepackage[T1]{fontenc}
\usepackage{lmodern}
\usepackage{amsmath,amssymb,amsthm,mathtools}
\usepackage{booktabs}
\usepackage{array}
\usepackage{enumitem}
\usepackage{microtype}
\usepackage[dvipsnames]{xcolor}
\usepackage[colorlinks=true,linkcolor=MidnightBlue,citecolor=MidnightBlue,urlcolor=MidnightBlue]{hyperref}

\numberwithin{equation}{section}
\newtheorem{theorem}{Theorem}[section]
\newtheorem{proposition}[theorem]{Proposition}
\newtheorem{lemma}[theorem]{Lemma}
\newtheorem{corollary}[theorem]{Corollary}
\theoremstyle{remark}
\newtheorem{remark}[theorem]{Remark}

\newcommand{\Z}{\mathbb Z}
\newcommand{\eps}{\varepsilon}
\newcommand{\Np}{N_{+}}
\newcommand{\Nm}{N_{-}}
\newcommand{\sq}{\square}
\setlist{itemsep=2pt,topsep=4pt}
\title[Gauss's mean value formulas for class numbers]{On the error term in Gauss's mean value formulas\\ for class numbers}

\author{Bingxin Jia}
\address{school of mathematical sciences\\ hebei normal university\\
	shijiazhuang, hebei 050024\\ china.}
\email{jiabx1010@163.com}

\author{CAIHUA LUO}
\address{school of science and engineering \\ the chinese university of hong kong \\ shenzhen, guangdong 518172\\ china.}
\email{luocaihua@cuhk.edu.cn}

\author{Xu Song}
\address{school of mathematical sciences\\ hebei normal university\\
	shijiazhuang, hebei 050024\\ china.}
\email{songx@hebtu.edu.cn}

\author{Wei Xiong}
\address{school of mathematics\\ hunan university\\ changsha, hunan 410082\\ China.}
\email{weixiong@amss.ac.cn}


\begin{document}

\begin{abstract}
Combining the frameworks of Chamizo--Iwaniec and Chamizo--Ubis with the recent character-sum and exponential-sum estimates of Bellotti and Kerr, we improve the error term in Gauss's mean value formulas for class numbers of binary quadratic forms from $O(N^{21/32+\eps})$ to $O(N^{55/84+\eps})$. If no Dirichlet $L$-function vanishes in $\mathrm{Re}\,s>7/8$, the exponent becomes $2274811/3487744$; under the Generalized Riemann Hypothesis for quadratic Dirichlet $L$-functions, it becomes $118/181$.
\end{abstract}

\maketitle

\tableofcontents

\section{Introduction}\label{sec:intro}

\subsection{Gauss's mean value formulas}
For an integer $D\equiv0,1\pmod 4$ that is not a perfect square, let $h(D)$ denote the number of $\mathrm{SL}_2(\Z)$-classes of primitive binary quadratic forms $ax^2+bxy+cy^2$ of discriminant $b^2-4ac=D$, positive definite when $D<0$; put $h(D)=0$ if $D\not\equiv0,1\pmod 4$. For positive non-square $D$ let $\varepsilon_D=(t+u\sqrt D)/2$, where $(t,u)$ is the least solution in positive integers of $t^2-Du^2=4$. Thus forms are counted up to proper equivalence, and $\varepsilon_D$ has norm $1$; for fundamental $D>0$, $h(D)$ is the narrow class number (see Section~\ref{ssec:classnumbers}). For positive square $D$ we set $h(D)\log\varepsilon_D=0$, so that square discriminants never contribute.

In Article~302 of the \emph{Disquisitiones Arithmeticae} \cite{Ga}, Gauss stated without proof an approximate formula for the average number of classes of forms of a given negative determinant. In modern notation his two-term formula reads
\begin{equation}\label{eq:gauss-neg}
\sum_{n\le N}h(-4n)=\frac{4\pi}{21\,\zeta(3)}N^{3/2}-\frac{2}{\pi^2}N+o(N).
\end{equation}
The first published proof is due to Lipschitz \cite{Li}. Vinogradov later obtained power savings in the error term by means of exponential sums, with bounds that also apply to the sphere problem \cite{Vi}. Chamizo and Iwaniec \cite{CI2}, using the methods they had developed for the sphere problem in \cite{CI1}, proved \eqref{eq:gauss-neg} with error term $O(N^{29/44+\eps})$. In a note added in proof they observed that Heath-Brown's estimates for the sphere problem \cite{HB} improve this to $O(N^{21/32+\eps})$, and this has remained the best known bound (see \cite[(1)]{CU}).

In Article~304, Gauss turned to positive determinants. He observed that the class number weighted by the logarithm of the fundamental unit behaves regularly, and proposed a formula of the shape $m\sqrt D-n$ for its average. Siegel \cite{Si} proved that
\[
\sum_{n\le N}h(4n)\log\varepsilon_{4n}=\frac{4\pi^2}{21\,\zeta(3)}N^{3/2}+O(N\log N).
\]
Shintani \cite{Sh} identified the secondary term, which is of order $N\log N$, with an error term $O(N^{\alpha})$ for every $\alpha>3/4$, and Chamizo and Ubis \cite{CU} reduced the error term to $O(N^{\alpha})$ for every $\alpha>21/32$.

In both cases $21/32$ is half of Heath-Brown's exponent $21/16$ for the sphere problem \cite{HB}: both proofs rest on the same two estimates as the sphere problem, namely a bound for sums of $L(1,\chi)$ over short intervals of discriminants, and a bound for exponential sums attached to ternary lattices. Recently Bellotti and Kerr \cite{BK} have improved the estimate for the exponent of the sphere problem from $21/16$ to $55/42$. The purpose of this paper is to transfer their improvement to Gauss's mean values.

\subsection{Main Results}
Now we state the main results of this paper.

\begin{theorem}\label{thm:neg}
For every $\eps>0$,
\begin{align*}
\sum_{n\le N}h(-4n)&=\frac{4\pi}{21\,\zeta(3)}\,N^{3/2}-\frac{2}{\pi^2}\,N+O_\eps\big(N^{55/84+\eps}\big),\\
\sum_{n\le N}h(-n)&=\frac{\pi}{18\,\zeta(3)}\,N^{3/2}-\frac{3}{2\pi^2}\,N+O_\eps\big(N^{55/84+\eps}\big).
\end{align*}
\end{theorem}

\begin{theorem}\label{thm:pos}
Let $C=\log(2\pi)-\zeta'(2)/\zeta(2)-1$ and $C'=C+\frac83\log 2$. For every $\eps>0$,
\begin{align*}
\sum_{n\le N}h(n)\log\varepsilon_n&=\frac{\pi^2}{18\,\zeta(3)}\,N^{3/2}-\frac{3}{\pi^2}\,(C+\log N)\,N+O_\eps\big(N^{55/84+\eps}\big),\\
\sum_{n\le N}h(4n)\log\varepsilon_{4n}&=\frac{4\pi^2}{21\,\zeta(3)}\,N^{3/2}-\frac{4}{\pi^2}\,(C'+\log N)\,N+O_\eps\big(N^{55/84+\eps}\big).
\end{align*}
\end{theorem}

The main and secondary terms are those of \cite{CI2} and \cite{CU}; the improvement is in the error exponent. The exponent $55/84=0.65476\ldots$ improves on $21/32=0.65625$ by $1/672$, half of the improvement $21/16-55/42=1/336$ obtained in \cite{BK} for the sphere problem, as one expects since $N=R^2$. The sphere estimate alone controls only one linear combination of sums of Hurwitz class numbers; our contribution is to establish the required estimates for the four sums separately, and Section~\ref{ssec:relation} describes the additional work.

Bellotti and Kerr \cite{BK} also showed that the non-vanishing of all Dirichlet $L$-functions  in the half-plane $\mathrm{Re}\,s>7/8$, announced in \cite{OAI}, would lower the sphere exponent to $2274811/1743872$. Under the same hypothesis we obtain the following.

\begin{theorem}\label{thm:cond}
Assume that no Dirichlet $L$-function has a zero in the half-plane $\mathrm{Re}\,s>7/8$. Then all four error terms in Theorems~\ref{thm:neg} and~\ref{thm:pos} are $O_\eps(N^{\vartheta+\eps})$ with
\[
\vartheta=\frac{2274811}{3487744}=0.6522299\ldots
\]
\end{theorem}

\begin{remark}\label{rem:grh}
If we assume the uniform generalized Lindel\"{o}f bound
\[
L(\tfrac12+it,\chi)\ll_\eps\big(q(1+|t|)\big)^\eps
\]
for primitive non-principal quadratic characters $\chi$ of conductor $q$, uniformly in $\chi$ and $t\in\mathbb R$ (this  follows from the Generalized Riemann Hypothesis for quadratic Dirichlet $L$-functions), then the four formulas of Theorems~\ref{thm:neg} and~\ref{thm:pos} hold, with the same main terms, with error terms $O_\eps(N^{118/181+\eps})$, where $118/181=0.6519337\ldots$ (see Section~\ref{ssec:lindelof}). 
\end{remark}

\subsection{Relation with the sphere problem}\label{ssec:relation}
Now we explain the relation of our results with  the Gauss sphere problem. Note that the sphere bound of \cite{BK} alone does not  imply Theorems~\ref{thm:neg} and~\ref{thm:pos}. Let $r_3(n)$ be the number of representations of $n$ as a sum of three squares. With the Hurwitz class number $\Nm$ of \eqref{eq:Nm}, Gauss's theorem on sums of three squares gives (see, e.g., \cite{Gro})
\[
r_3(n)=12\big(\Nm(4n)-2\Nm(n)\big).
\]
Thus the sphere problem controls only one linear combination of $\sum_{n\le x}\Nm(4n)$ and $\sum_{n\le x}\Nm(n)$, whose main terms in Proposition~\ref{prop:Nsums} combine to $\frac{4\pi}3x^{3/2}$, with no linear term. Theorem~\ref{thm:neg} requires the two sums separately, and the regulator-weighted class numbers $\Np(n)$ of Theorem~\ref{thm:pos} do not occur in this identity at all.

We therefore establish the estimates underlying the sphere bound for the four sequences $\Np(n)$, $\Np(4n)$, $\Nm(n)$ and $\Nm(4n)$, within the frameworks of Chamizo and Iwaniec \cite{CI2} for negative discriminants and of Chamizo and Ubis \cite{CU} for positive ones. The sharper character-sum and exponential-sum estimates are due to Bellotti and Kerr \cite{BK}.

On the arithmetic side, we treat both signs of the discriminant and the residue classes modulo $8$ required by the four sequences (see Proposition~\ref{prop:boundary}). For positive discriminants, square discriminants are excluded from the class-number sums, but square terms reappear once the $L$-values are truncated to finite character sums, and they give principal characters in the character sums of Section~\ref{sec:charsums}; their total contribution is $O(R^{1+\eps})$.

On the analytic side, Gauss's reduced forms and Siegel's geodesic weights reduce the four dual sums to double exponential sums (see Lemma~\ref{lem:reduction}). Bellotti and Kerr state that their estimates hold for subintervals \cite[Section 4.2]{BK}; Lemma~\ref{lem:BK41} gives an explicit completion argument for the shift-dependent intervals arising here. Shintani's summation formulas, in the form of \cite{CU} for positive discriminants and derived here for negative ones (see Lemma~\ref{lem:sum-neg}), and M\"{o}bius inversion then give the four mean-value formulas with their known main and secondary terms.

\subsection{Outline of the proof}
After M\"{o}bius inversion it suffices to estimate the sums of $N_\pm(n)$ and $N_\pm(4n)$. We smooth the cut-off at $R=\sqrt N$ over an interval of width $H$, and the proof separates into two estimates. The boundary contribution is a short sum of values $L(1,\chi)$; it is bounded in Section~\ref{sec:boundary} by the character-sum estimates of Section~\ref{sec:charsums}, which come from \cite[Section 3]{BK} and rest on Heath-Brown's quadratic large sieve \cite{HB95} and the Weyl bound of Petrow and Young \cite{PY}. Shintani's summation formulas express the smoothed sums as main terms plus dual sums (see Section~\ref{sec:smoothed}), and the dual sums are controlled by the exponential-sum estimates of Section~\ref{sec:expsums}, where Bourgain's exponent pair \cite{Bou} enters. Both estimates are stated for a general width $H$ and general character-sum and exponential-sum (see Remark~\ref{rem:general-char} and Lemma~\ref{lem:dual}). Proposition~\ref{prop:Nsums} and Corollary~\ref{cor:transfer} combine them and recover the sharp sums; this shows explicitly how the final error exponent depends on the character-sum and exponential-sum bounds used. For the unconditional result, $H=R^{-13/21}$ gives the error $O(R^{55/42+\eps})$, that is, $O(N^{55/84+\eps})$. The conditional results follow from the same transfer, with other ingredients and other choices of $H$ (see Table~\ref{tab:opt}). Section~\ref{sec:prelim} collects the class-number identities and summation formulas used throughout.

\subsection{Notation}
We write $e(x)=e^{2\pi ix}$. For an integer $D\equiv0,1\pmod4$, $\chi_D(m)=(D/m)$ denotes the Kronecker symbol; for odd $m$ and any integer $a$, $(a/m)$ denotes the Jacobi symbol. We write $m=\sq$ to mean that $m$ is a perfect square. For a positive integer $k$, $\mu(k)$ denotes the M\"{o}bius function.  The letter $\eps$ denotes an arbitrarily small positive number and $R$ denotes a large parameter.

\section{Class numbers and summation formulas}\label{sec:prelim}

In this section, we will collect the following  results: M\"{o}bius inversion and Dirichlet's class number formula (see Section~\ref{ssec:classnumbers}), the lattice point interpretations of the class numbers (see Section~\ref{ssec:lattice}), the summation formulas derived from Shintani's functional equation (see Section~\ref{ssec:shintani}), and the smooth weights (see Section~\ref{ssec:smoothing}).

\subsection{Class numbers}\label{ssec:classnumbers}
Following \cite{CU}, we work with class numbers summed over all forms, primitive or not. For $n\ge1$ put
\begin{equation}\label{eq:Nm}
\Np(n)=\sum_{k^2\mid n}h(n/k^2)\log\varepsilon_{n/k^2},\qquad
\Nm(n)=\sum_{k^2\mid n}\frac{h(-n/k^2)}{\delta_{-n/k^2}},
\end{equation}
where $\delta_{-3}=3$, $\delta_{-4}=2$ and $\delta_D=1$ for all other $D$. Thus $\Nm(n)$ is the Hurwitz class number of discriminant $-n$. The function $\Np$ is supported on $n\equiv0,1\pmod4$, and $\Nm$ on $n\equiv0,3\pmod4$. The class numbers in Theorems~\ref{thm:neg} and~\ref{thm:pos} are recovered by M\"{o}bius inversion.

\begin{lemma}\label{lem:mobius}
For $n\ge1$,
\begin{align*}
h(n)\log\varepsilon_n&=\sum_{k^2\mid n}\mu(k)\Np(n/k^2), & h(4n)\log\varepsilon_{4n}&=\sum_{\substack{k^2\mid n\\2\nmid k}}\mu(k)\big(\Np(4n/k^2)-\Np(n/k^2)\big),\\
\frac{h(-n)}{\delta_{-n}}&=\sum_{k^2\mid n}\mu(k)\Nm(n/k^2), & \frac{h(-4n)}{\delta_{-4n}}&=\sum_{\substack{k^2\mid n\\2\nmid k}}\mu(k)\big(\Nm(4n/k^2)-\Nm(n/k^2)\big).
\end{align*}
\end{lemma}

\begin{proof}
The first identity in each row is M\"{o}bius inversion of \eqref{eq:Nm} over square divisors. For the second, write $k=2k'$ for even $k$; since $\mu(2k')=-\mu(k')\mathbf 1_{2\nmid k'}$,
\[
h(4n)\log\varepsilon_{4n}=\sum_{k^2\mid4n}\mu(k)\Np\Big(\frac{4n}{k^2}\Big)=\sum_{\substack{k^2\mid n\\2\nmid k}}\mu(k)\Np\Big(\frac{4n}{k^2}\Big)-\sum_{\substack{k'^2\mid n\\2\nmid k'}}\mu(k')\Np\Big(\frac n{k'^2}\Big),
\]
and the same computation applies to $\Nm$ \cite[Lemma 2.1]{CU}.
\end{proof}

For every non-square discriminant $D$ (fundamental or not), we have Dirichlet's class number formula:
\[\begin{aligned}
h(D)\log\varepsilon_D=\sqrt D\,L(1,\chi_D)  \quad &\text{for $D>0$},\\
h(D)/\delta_D=\sqrt{|D|}\,L(1,\chi_D)/\pi   \quad &\text{ for $D<0$},
\end{aligned}\]
where $\chi_D$ is the Kronecker symbol and $L(1, \chi_D)$ is the value of the Dirichlet $L$-function $L(s, \chi_D)$ at $1$.
Summing over square divisors, we obtain
\begin{equation}\label{eq:Lform}
\Np(n)=\sum_{\substack{d^2\mid n\\ n/d^2\in\mathcal D_+}}\sqrt{n/d^2}\,L(1,\chi_{n/d^2}),\qquad
\Nm(n)=\frac1\pi\sum_{\substack{d^2\mid n\\ n/d^2\in\mathcal D_-}}\sqrt{n/d^2}\,L(1,\chi_{-n/d^2}),
\end{equation}
where $\mathcal D_+$ is the set of non-square $a\ge1$ with $a\equiv0,1\pmod 4$, and $\mathcal D_-$ the set of $a\ge1$ with $a\equiv0,3\pmod4$.

For a fundamental discriminant $D>0$, $h(D)$ is the narrow class number of $K=\mathbb Q(\sqrt D)$, and
\[
h(D)\log\varepsilon_D=2h_K\log u_K,
\]
where $h_K$ is the class number of $K$ and $u_K>1$ the fundamental unit of $K$.

\subsection{Lattice point representations}\label{ssec:lattice}

In Section~\ref{sec:expsums}, we will need two classical representations of $\Nm$ and $\Np$ as weighted counts of lattice points.
The first is Gauss's count of reduced forms, in the form of \cite[Lemma 2.1]{CI2} (see also \cite[p.~77]{CU}).

\begin{lemma}\label{lem:gauss}
For $n\ge1$, $\Nm(n)$ equals the number of $(a,b,c)\in\Z^3$ with $4ac-b^2=n$ and either $-a<b\le a<c$ or $0\le b\le a=c$, except that $\frac23$ must be subtracted if $n=3k^2$ and $\frac12$ if $n=4k^2$.
\end{lemma}

For an indefinite form $(a,b,c)$ with $a>0$ and non-square discriminant, let $\mu(a,b,c)$ be the hyperbolic length of the intersection of the geodesic $a(x^2+y^2)+bx+c=0$ in the upper half-plane with the standard fundamental domain $\{x+iy:\ |x|<\frac12,\ x^2+y^2>1\}$. The second representation goes back to Siegel \cite{Si}; we use it in the form of \cite[Lemmas 3.1, 3.2]{CU}.

\begin{lemma}\label{lem:siegel}
For non-square $n\ge1$,
\[
\Np(n)=\sum_{\substack{b^2-4ac=n\\ a>0}}\mu(a,b,c).
\]
If $\mu(a,b,c)\ne0$, then $a+c<|b|/2$, $a\le\sqrt{n/3}$, $|b|\le2\sqrt{n/3}$ and $a|c|\le n/4$. Moreover $\mu(a,b,c)\ll n^\eps$, and for fixed $a,c$ the set of $b$ with $\mu(a,b,c)\ne0$ is a union of $O(1)$ intervals, on each of which $\mu(a,\cdot,c)$ is monotone (see \cite[Proof of Proposition 3.3]{CU}).
\end{lemma}

\subsection{Shintani's functional equation and summation formulas}\label{ssec:shintani}
For positive discriminants we use the summation formulas of Chamizo and Ubis \cite[Propositions 2.3, 2.4]{CU}, which they derive from Shintani's functional equation; we recall them in Lemma~\ref{lem:sum-pos}. For $g\in C_0^\infty((0,\infty))$  put
\[
\tilde g(x)=\int_0^\infty g(t)\cos(\pi xt)\,dt,\qquad \breve g(x)=\int_0^\infty g(t)\sin(\pi xt)\,dt,\qquad G(s)=\int_0^\infty g(t)t^{s-1}\,dt .
\]

\begin{lemma}\label{lem:sum-pos}
For $g\in C_0^\infty((0,\infty))$,
\begin{align*}
\sum_{n\ge1}\frac{\Np(n)}{\sqrt n}g(\sqrt n)&=M_1(g)+2\sum_{n\ge1}\frac{\Np(4n)}{\sqrt{4n}}\,\tilde g(\sqrt{4n})+E_1(g),\\
\sum_{n\ge1}\frac{\Np(4n)}{\sqrt{4n}}g(\sqrt{4n})&=M_4(g)+\sum_{n\ge1}\frac{\Np(n)}{\sqrt{n}}\,\tilde g(\sqrt{n})+E_4(g),
\end{align*}
where
\begin{align*}
M_1(g)&=\frac{\pi^2}6\int_0^\infty tg(t)\,dt-\int_0^\infty g(t)\log(2\pi t)\,dt-\sum_{n\ge1}g(n)\log n+\sum_{d\ge1}\frac{\Lambda(d)}d\sum_{n\ge1}g(dn),\\
M_4(g)&=\frac{\pi^2}{12}\int_0^\infty tg(t)\,dt-\frac12\int_0^\infty g(t)\log(2\pi t)\,dt-\sum_{n\ge1}g(2n)\log n+\sum_{d\ge1}\frac{\Lambda(d)}d\sum_{n\ge1}g(2dn)\\
&\qquad-\log2\sum_{k\ge1}2^{-k}\sum_{n\ge1}g(2^kn),
\end{align*}
and
\begin{align*}
E_1(g)&=2\sum_{n\ge1}\beta_n\tilde g(n)+2\pi\sum_{n\ge1}\frac{\Nm(4n)}{\sqrt{4n}}\int_0^\infty g(t)e^{-\pi\sqrt{4n}\,t}\,dt+2\int_0^1\sum_{n\ge1}\tilde g(2n/t)\frac{dt}{t(1+t)},\\
E_4(g)&=\sum_{n\ge1}\alpha_n\tilde g(n)+\pi\sum_{n\ge1}\frac{\Nm(n)}{\sqrt{n}}\int_0^\infty g(t)e^{-\pi\sqrt{n}\,t}\,dt-\int_0^1\sum_{n\ge1}\tilde g(n/t)\frac{dt}{t(1+t)}.
\end{align*}
Here $\sum_n\alpha_nn^{-s}=\zeta(s)\big(\zeta'(s+1)/\zeta(s+1)-\zeta'(s)/\zeta(s)\big)$ and $\sum_n\beta_nn^{-s}=\sum_n\alpha_n(2n)^{-s}+\zeta(s)\log2/(2^{s+1}-1)$; in particular $\alpha_n,\beta_n\ll\log(2n)$.
\end{lemma}

For negative discriminants we need analogous formulas. They follow from the second coordinate of Shintani's vector functional equation, which, unlike the first, involves only negative discriminants. In the notation of \cite[Section 2]{CU}, put
\[
\eta_1(s)=\sum_{n\ge1}\Nm(4n)(4n)^{-(s+1)/2},\qquad \eta_2(s)=\sum_{n\ge1}\Nm(n)\,n^{-(s+1)/2}.
\]
These series converge absolutely for $\mathrm{Re}\,s>2$. By \cite[Theorem 2.2]{CU} (which is \cite[Theorem 2]{Sh} after a change of variables), $\eta_1$ and $\eta_2$ extend to meromorphic functions of order one. Their only poles are simple poles at $s=2$ and $s=1$, with residues $\pi/12$ and $-1/4$ for $\eta_1$, and $\pi/6$ and $-1/2$ for $\eta_2$. Moreover,
\begin{equation}\label{eq:FE}
\eta_2(1-s)=-2\pi^{-s}\Gamma(s)\sin\Big(\frac{\pi s}2\Big)\eta_1(s)-(2\pi)^{-s}\Gamma(s)\zeta(s).
\end{equation}
Indeed, the matrix $A(s)$ of \cite[Theorem 2.2]{CU} has a zero lower-left entry, and the second coordinate of $\vec B(s)$ there is $\sec(\pi s/2)$, so the second coordinate of the vector equation involves only $\eta_1$, $\eta_2$ and $\zeta$. Replacing $s$ by $1-s$ in \eqref{eq:FE} and using $\Gamma(s)\Gamma(1-s)=\pi/\sin(\pi s)$ and the functional equation of $\zeta$ gives
\begin{equation}\label{eq:FE2}
\eta_1(1-s)=-\pi^{-s}\Gamma(s)\sin\Big(\frac{\pi s}2\Big)\eta_2(s)-\frac12\pi^{-s}\Gamma(s)\zeta(s).
\end{equation}

\begin{lemma}\label{lem:sum-neg}
For $g\in C_0^\infty((0,\infty))$,
\begin{align*}
\sum_{n\ge1}\frac{\Nm(n)}{\sqrt n}g(\sqrt n)&=\frac\pi6\int_0^\infty tg(t)\,dt-\frac12\int_0^\infty g(t)\,dt-2\sum_{n\ge1}\frac{\Nm(4n)}{\sqrt{4n}}\,\breve g(\sqrt{4n})-\int_0^\infty\frac{g(t)\,dt}{e^{2\pi t}-1},\\
\sum_{n\ge1}\frac{\Nm(4n)}{\sqrt{4n}}g(\sqrt{4n})&=\frac\pi{12}\int_0^\infty tg(t)\,dt-\frac14\int_0^\infty g(t)\,dt-\sum_{n\ge1}\frac{\Nm(n)}{\sqrt n}\,\breve g(\sqrt n)-\frac12\int_0^\infty\frac{g(t)\,dt}{e^{\pi t}-1}.
\end{align*}
\end{lemma}

\begin{proof}
We argue as in \cite[Proof of Proposition 2.3]{CU}. Here $G$ is entire and $G(s)\ll_A(1+|s|)^{-A}$ in vertical strips, and by \eqref{eq:FE} and convexity $\eta_2$ grows at most polynomially in $|\mathrm{Im}\,s|$ for $|\mathrm{Re}\,s|\le3$. For $2<c<3$ and $-2<c'<-1$,
\[
\sum_{n\ge1}\frac{\Nm(n)}{\sqrt n}g(\sqrt n)=\frac1{2\pi i}\int_{(c)}G(s)\eta_2(s)\,ds=\frac\pi6G(2)-\frac12G(1)+\frac1{2\pi i}\int_{(1-c')}G(1-s)\eta_2(1-s)\,ds,
\]
with $G(2)=\int tg$ and $G(1)=\int g$. On $\mathrm{Re}\,s=1-c'>2$ insert \eqref{eq:FE} and expand
\[
\eta_1(s)=\sum_{n\ge1}\frac{\Nm(4n)}{\sqrt{4n}}(4n)^{-s/2},\qquad\zeta(s)=\sum_{m\ge1}m^{-s}.
\]
Moving each term to $0<\mathrm{Re}\,s<1$ and using
\[
\Gamma(s)\sin\frac{\pi s}2=\int_0^\infty t^{s-1}\sin t\,dt,\qquad\Gamma(s)(2\pi m)^{-s}=\int_0^\infty t^{s-1}e^{-2\pi mt}\,dt,
\]
Mellin inversion gives
\begin{align*}
\frac1{2\pi i}\int G(1-s)\big(\pi\sqrt{4n}\big)^{-s}\Gamma(s)\sin\frac{\pi s}2\,ds&=\breve g(\sqrt{4n}),\\
\sum_{m\ge1}\frac1{2\pi i}\int G(1-s)\Gamma(s)(2\pi m)^{-s}\,ds&=\int_0^\infty\frac{g(t)\,dt}{e^{2\pi t}-1}.
\end{align*}
Substituting these identities into the shifted contour integral gives the first formula. The second formula follows from \eqref{eq:FE2} similarly.
\end{proof}

\begin{remark}
The residues of $\eta_1,\eta_2$ determine the main terms of our theorems; they are given by \cite[Theorem 2.2]{CU}. For the error terms only the boundedness of the constants in front of the dual sums matters.
\end{remark}

\subsection{Smooth weights}\label{ssec:smoothing}
Choose a nonincreasing $\rho_0\in C^\infty(\mathbb R)$ with $\rho_0=1$ on $(-\infty,0]$ and $\rho_0=0$ on $[1,\infty)$. Fix also $\kappa\in C^\infty([0,\infty))$ with $0\le\kappa\le1$, $\kappa=0$ on $[0,\frac14]$ and $\kappa=1$ on $[\frac12,\infty)$. For $R\ge16$ and $0<H\le\frac14$ put
\[
f(t)=\rho_0\Big(\frac{t-R}H\Big),\qquad g(t)=t\,f(t)\,\kappa(t).
\]
Then $f=1$ on $[0,R]$, $f=0$ on $[R+H,\infty)$, and $f^{(j)}\ll_jH^{-j}$. Moreover $g\in C_0^\infty((0,\infty))$, $0\le g(t)\le t$, $g=0$ on $[0,\frac14]$, $g(t)=t$ on $[\frac12,R]$, and $g(\sqrt n)/\sqrt n=f(\sqrt n)$ for every $n\ge1$. (For $H=R^{-13/21}$ or $H=R^{-3/5}$, the condition $H\le\frac14$ holds since $R\ge16$; for $H=R^{-71/181}$ it holds when $R\ge64$.) All derivatives of $g$ of order at least two are supported in $[\frac14,\frac12]\cup[R,R+H]$; on $[\frac14,\frac12]$ they do not depend on $R$ or $H$.

\begin{lemma}\label{lem:transforms}
Let $F(t)=tf(t)$, and for $\xi\ge\frac12$ put $I(\xi)=\int_0^\infty g(t)\sin(2\pi\xi t)\,dt$ and $J(\xi)=\int_0^\infty g(t)\cos(2\pi\xi t)\,dt$. Then, for every $A>0$,
\begin{align*}
I(\xi)&=\frac1{2\pi\xi}\int_R^{R+H}tf'(t)\cos(2\pi\xi t)\,dt-\frac1{(2\pi\xi)^2}\int_R^{R+H}f'(t)\sin(2\pi\xi t)\,dt+O_A(\xi^{-A}),\\
J(\xi)&=-\frac1{2\pi\xi}\int_R^{R+H}tf'(t)\sin(2\pi\xi t)\,dt-\frac1{(2\pi\xi)^2}\int_R^{R+H}f'(t)\cos(2\pi\xi t)\,dt+O_A(\xi^{-A}),
\end{align*}
and
\[
\begin{Bmatrix}I(\xi)\\ J(\xi)\end{Bmatrix}=-\frac1{(2\pi\xi)^2}\int_R^{R+H}F''(t)\begin{Bmatrix}\sin(2\pi\xi t)\\ \cos(2\pi\xi t)\end{Bmatrix}dt+O_A(\xi^{-A}).
\]
In particular $I(\xi),J(\xi)\ll_A\xi^{-A}+R\min\{\xi^{-1},H^{-1}\xi^{-2}\}$, and also $I(\xi),J(\xi)\ll RH^{-3}\xi^{-4}$.
\end{lemma}

\begin{proof}
We treat $J$; the case of $I$ is analogous, and its first formula is \cite[(4.16)]{BK}. Since $g'=(t\kappa)'$ on $[0,R]$ and $g'=f+tf'$ on $[R,R+H]$,
\[
J(\xi)=-\frac1{2\pi\xi}\int_0^\infty g'(t)\sin(2\pi\xi t)\,dt,
\]
and
\begin{align*}
\int_0^R(t\kappa)'\sin(2\pi\xi t)\,dt&=-\frac{\cos(2\pi\xi R)}{2\pi\xi}+\frac1{2\pi\xi}\int_0^R(t\kappa)''\cos(2\pi\xi t)\,dt,\\
\int_R^{R+H}f\sin(2\pi\xi t)\,dt&=\frac{\cos(2\pi\xi R)}{2\pi\xi}+\frac1{2\pi\xi}\int_R^{R+H}f'\cos(2\pi\xi t)\,dt.
\end{align*}
The boundary terms cancel, and $(t\kappa)''\in C_0^\infty([\frac14,\frac12])$ does not depend on $R,H$, so its Fourier transform is $O_A(\xi^{-A})$. This gives the first two formulas; the third follows by integrating by parts twice, since $g''=(t\kappa)''$ on $[\frac14,\frac12]$ and $g''=F''$ on $[R,R+H]$. The bounds follow from
\[
\|tf'\|_1\ll R,\qquad\|f'\|_1\ll1,\qquad\|F''\|_1\ll\frac RH,\qquad\|F^{(4)}\|_1\ll\frac R{H^3},
\]
the last after two further integrations by parts.
\end{proof}

\section{Character sums over short intervals}\label{sec:charsums}

The boundary terms of Section~\ref{sec:boundary} lead to the character sums defined below. We extend the sums of \cite[Section 3.3]{BK} to both signs of the discriminant. As in \cite[(3.13)]{BK}, square values of $n$ are included; for the positive sign they give principal characters, and they are treated separately in Lemma~\ref{lem:PY}. We first record the two estimates of \cite[Section 3]{BK} in this generality. Let $\Phi$ be smooth and supported in $[\frac23,\frac53]$, and let $W$ be smooth and supported in $[\frac12,2]$, with $\Phi^{(r)},W^{(r)}\ll_r1$ uniformly. All estimates below are uniform in such weights. For $\sigma\in\{1,-1\}$, a residue class $\nu\bmod8$, $X>0$ and $L,J\ge1$, put
\begin{equation}\label{eq:S}
S^\sigma_\nu(X;L,J)=\frac1J\sum_{\substack{m\ \mathrm{odd}\\ m\ne\sq}}W\Big(\frac mJ\Big)\sum_{n\equiv\nu\ (8)}\Phi\Big(\frac{n-X}L\Big)\Big(\frac{\sigma n}m\Big).
\end{equation}
For odd $m$ we have $(-n/m)=(-4n/m)$, and $(-4n/m)=0$ for even $m$. Hence $S^{-1}_\nu$ is exactly the sum $S_\nu$ of \cite[(3.13)]{BK}.

\begin{lemma}\label{lem:BKprop}
If $X,L,J\le R^{O(1)}$, then $S^\sigma_\nu(X;L,J)\ll_\eps R^\eps\big((JL)^{1/4}+J^{1/2}L^{-1/4}\big)$.
\end{lemma}

\begin{proof}
For $\sigma=-1$ this is \cite[Proposition 3.2]{BK}. Let $\sigma=1$ and $m=bc^2$ with $b>1$ odd and squarefree. Then
\[
\Big(\frac n{bc^2}\Big)=\Big(\frac nb\Big)\mathbf 1_{(n,c)=1}=\sum_{e\mid(n,c)}\mu(e)\Big(\frac eb\Big)\Big(\frac{n/e}b\Big),
\]
so \cite[(3.19)]{BK} holds with $(-4e/b)$ replaced by $(e/b)$. The argument there uses only $|(e/b)|\le1$: the sum over $b$ is bounded by the triangle inequality, and \cite[Lemma 3.1]{BK} is applied.
\end{proof}

\begin{lemma}\label{lem:PY}
Let $1\le Q\le R$, $2\le L\le Q$ and $1\le J\le R^{O(1)}$. Then
\[
S^\sigma_\nu(Q^2;L,J)\ll_\eps R^\eps Q^{1/3}LJ^{-1/2}+\mathbf 1_{\sigma=1}.
\]
\end{lemma}

\begin{proof}
The $n$-sum has $O(L)$ terms, all with $n\asymp Q^2$. If $\sigma n\ne\sq$, then $(\sigma n/m)=\chi_{4\sigma n}(m)$ for odd $m$, $\chi_{4\sigma n}(m)=0$ for even $m$, and $\chi_{4\sigma n}$ is non-principal, induced by a primitive character of conductor at most $4n$. By \cite[Lemma 3.3]{BK}, whose proof rests on the Weyl bound \cite[Theorem 1.1]{PY},
\[
\sum_{\substack{m\ \mathrm{odd}\\ m\ne\sq}}W\Big(\frac mJ\Big)\Big(\frac{\sigma n}m\Big)=\sum_{m\ne\sq}\chi_{4\sigma n}(m)W\Big(\frac mJ\Big)\ll_\eps J^{1/2}n^{1/6+\eps}.
\]
If $\sigma=1$ and $n=\sq$, the $m$-sum is $O(J)$. Consecutive squares exceeding $Q^2$ differ by more than $2Q$, while $\Phi((n-Q^2)/L)$ is supported on an interval of length $L\le Q$; so there is at most one such $n$. Hence
\[
S^\sigma_\nu(Q^2;L,J)\ll_\eps\frac1J\big(L\,J^{1/2}Q^{1/3+\eps}+\mathbf 1_{\sigma=1}J\big).
\]
\end{proof}

\begin{lemma}\label{lem:combine}
Under the hypotheses of Lemma~\ref{lem:PY},
\[
S^\sigma_\nu(Q^2;L,J)\ll_\eps R^\eps\big(Q^{1/9}L^{1/2}+Q^{1/6}L^{3/8}+1\big)
\]
and
\[
S^\sigma_\nu(Q^2;L,J)\ll_\eps R^\eps\big((JL)^{1/4}+Q^{1/15}J^{3/10}+1\big).
\]
\end{lemma}

\begin{proof}
By Lemmas~\ref{lem:BKprop} and~\ref{lem:PY}, $S^\sigma_\nu(Q^2;L,J)\ll_\eps R^\eps(\min\{A+B,C\}+1)$, where
\[
A=(JL)^{1/4},\qquad B=J^{1/2}L^{-1/4},\qquad C=Q^{1/3}LJ^{-1/2}.
\]
For $A,B,C\ge0$, as in \cite[(3.23)]{BK} and \cite[(6.6)]{BK},
\[
\min\{A+B,C\}\le A^{2/3}C^{1/3}+B^{1/2}C^{1/2},\qquad\min\{A+B,C\}\le A+B^{4/5}C^{1/5},
\]
and
\[
A^{2/3}C^{1/3}=Q^{1/9}L^{1/2},\qquad B^{1/2}C^{1/2}=Q^{1/6}L^{3/8},\qquad B^{4/5}C^{1/5}=Q^{1/15}J^{3/10}.\qedhere
\]
\end{proof}

\section{Boundary terms: short sums of \texorpdfstring{$L(1,\chi)$}{L(1,chi)}}\label{sec:boundary}

In this section we evaluate the boundary terms, which are short sums of $L(1,\chi)$ over discriminants. To treat the four cases needed for Theorems~\ref{thm:neg} and~\ref{thm:pos} at once, we work in the following general setting. Let $\sigma\in\{1,-1\}$. For each $d\ge1$ let $\mathcal A_d$ be a set of residue classes mod $8$ that depends only on the parity of $d$, such that $\sigma a\equiv0,1\pmod4$ for every $a$ in these classes. We call the family $(\mathcal A_d)$ \emph{balanced} if, for each $d$, the classes of $\mathcal A_d$ with $\sigma a\equiv1\pmod 8$ and those with $\sigma a\equiv5\pmod8$ have the same cardinality. For $b\ge1$ let $\delta_d(b)$ be the density of the set of $a\ge1$ with $a\bmod8\in\mathcal A_d$ and $(a,b)=1$, and put
\begin{equation}\label{eq:cconst}
c(\sigma,\mathcal A)=\sum_{d\ge1}\frac2{d^3}\sum_{b\ge1}\frac{\delta_d(b)}{b^2}.
\end{equation}
For $R\ge16$ and $0<H\le\frac14$ define
\[
\mathcal B(R)=\sum_{d\ge1}\ \sum_{\substack{R^2/d^2<a\le(R+H)^2/d^2\\ a\bmod 8\in\mathcal A_d,\ \sigma a\ne\sq}}\sqrt a\,f(d\sqrt a)\,L(1,\chi_{\sigma a}).
\]

\begin{proposition}\label{prop:boundary}
If $(\mathcal A_d)$ is balanced, then
\[
\mathcal B(R)=c(\sigma,\mathcal A)\int_R^{R+H}t^2f(t)\,dt+O_\eps\Big(R^\eps\big(R+R^{10/9}(RH)^{1/2}+R^{7/6}(RH)^{3/8}\big)\Big).
\]
\end{proposition}

\begin{proof}
We follow \cite[Section 3]{BK}. Let $Z_d=16R^2/d^2$, and let $U\in C^\infty([0,\infty))$ with $U=1$ on $[0,1]$ and $U=0$ on $[2,\infty)$. Every $\chi_{\sigma a}$ in $\mathcal B(R)$ is non-principal of modulus $|\sigma a|\le Z_d$, and for each $d\le R+H$ there are $O(1+RH/d^2)$ values of $a$. As in \cite[(3.3)--(3.4)]{BK},
\[
\sqrt a\,L(1,\chi_{\sigma a})=\sqrt a\sum_{m\ge1}\frac{\chi_{\sigma a}(m)}mU\Big(\frac m{Z_d}\Big)+O(\log R),
\]
so that, with $a$ running as in $\mathcal B(R)$,
\begin{equation}\label{eq:B-trunc}
\mathcal B(R)=\sum_{d\ge1}\sum_{m\ge1}\frac1mU\Big(\frac m{Z_d}\Big)\sum_a\sqrt a\,f(d\sqrt a)\,\chi_{\sigma a}(m)+O(R\log R).
\end{equation}
We split $m$ into $b^2$, $2b^2$, and $2^km_0$ with $k\ge0$, $m_0$ odd and $m_0\ne\sq$. The moduli $b^2$ and $2b^2$ give the main term \eqref{eq:B-main}; the moduli $2^km_0$ give the error term \eqref{eq:B-rest}.

For $D=\sigma a$,
\[
\chi_D(b^2)=\mathbf 1_{(a,b)=1},\qquad\chi_D(2b^2)=\chi_D(2)\mathbf 1_{(a,b)=1},\qquad\chi_D(2)=\begin{cases}0,&D\equiv0\pmod4,\\ \pm1,&D\equiv1,5\pmod8,\end{cases}
\]
so $\chi_D(2)$ depends only on $a\bmod8$. By \cite[(3.5)]{BK}, $U(b^2/Z_d)$ may be replaced by $1$ at a cost $O(R)$. Summing the smooth function $x\mapsto\sqrt xf(d\sqrt x)$ over residue classes, as in \cite[(3.8)]{BK},
\begin{equation}\label{eq:B-sq}
\sum_a\sqrt a\,f(d\sqrt a)\mathbf 1_{(a,b)=1}=\frac{2\delta_d(b)}{d^3}\int_R^{R+H}t^2f(t)\,dt+O\Big(\frac{\tau(b)R}d\Big),
\end{equation}
since $\int_{R^2/d^2}^{(R+H)^2/d^2}\sqrt xf(d\sqrt x)\,dx=2d^{-3}\int_R^{R+H}t^2f(t)\,dt$. The sum with the extra factor $\chi_{\sigma a}(2)$ is $O(\tau(b)R/d)$, because $(\mathcal A_d)$ is balanced and the classes $a\bmod8$ coprime to $b$ have equal densities. If $\sigma=1$, the window $(R^2/d^2,(R+H)^2/d^2]$ has length $\le3RH/d^2<R/d$, less than the gap between consecutive squares there; so it contains at most one square, and excluding it changes \eqref{eq:B-sq} by $O(R/d)$. Summing over $b$ with the weights $b^{-2}$ and over $d$, the terms $m=b^2,2b^2$ of \eqref{eq:B-trunc} equal
\begin{equation}\label{eq:B-main}
c(\sigma,\mathcal A)\int_R^{R+H}t^2f(t)\,dt+O(R\log R).
\end{equation}

Let $\chi_\nu(2)$ denote the common value of $\chi_{\sigma a}(2)$ for $a\equiv\nu\pmod8$. For $m=2^km_0$ and $a\equiv\nu\pmod8$,
\[
\chi_{\sigma a}(m)=\chi_\nu(2)^k\Big(\frac{\sigma a}{m_0}\Big).
\]
If $\sigma=1$, we add back the square values of $a$ (at most one for each $d$), at a cost
\[
\ll\sum_{d\le R+H}\frac Rd\sum_{m\le2Z_d}\frac1m\ll R\log^2R.
\]
Insert the partitions of unity of \cite[Section 3.3]{BK}: $\vartheta$ and dyadic $\Psi(x/L)$ in $x=a-R^2/d^2$, and dyadic $w(m/J)$. The $\vartheta$-term has $O(1)$ values of $a$ for each $d$ and contributes $O(R\log^2R)$. With $Q=R/d$, $J_k=J/2^k$, and $\Phi,W$ as in \cite[Section 3.3]{BK} (they satisfy the hypotheses of Section~\ref{sec:charsums}), the contribution of $(d,\nu,k,L,J)$ is, as in \cite[(3.14)]{BK},
\begin{equation}\label{eq:B-piece}
\chi_\nu(2)^k\,Q\,2^{-k}\,S^\sigma_\nu(Q^2;L,J_k),\qquad2\le L\le\frac{4RH}{d^2}\le Q,\quad J\ll Z_d,
\end{equation}
and it vanishes if $J_k<1$, since $m_0=1$ is a square. By Lemma~\ref{lem:combine}, \eqref{eq:B-piece} is
\[
\ll R^\eps2^{-k}\big(Q^{10/9}L^{1/2}+Q^{7/6}L^{3/8}+Q\big).
\]
Summing over $k\ge0$, over the $O(\log R)$ dyadic $J$, over dyadic $L\le4RH/d^2$, and over $d$, the terms $m=2^km_0$ contribute
\begin{equation}\label{eq:B-rest}
\ll R^\eps\Big(R+R^{10/9}(RH)^{1/2}\sum_{d\ge1}d^{-19/9}+R^{7/6}(RH)^{3/8}\sum_{d\ge1}d^{-23/12}\Big).
\end{equation}
The proposition follows from \eqref{eq:B-trunc}, \eqref{eq:B-main} and \eqref{eq:B-rest}.
\end{proof}

\begin{remark}\label{rem:general-char}
The proof uses Lemma~\ref{lem:combine} only to bound \eqref{eq:B-piece}. Suppose instead that, under the hypotheses of Lemma~\ref{lem:PY}, $S^\sigma_\nu(Q^2;L,J)\ll_\eps R^\eps\big(1+Q^{\beta}L^{\gamma}\big)$ for some $\beta,\gamma\ge0$ with $\beta+2\gamma>0$. Then the same argument gives Proposition~\ref{prop:boundary} with error term $O_\eps\big(R^\eps(R+R^{1+\beta}(RH)^{\gamma})\big)$, since the sum over $d$ is now $\sum_dd^{-1-\beta-2\gamma}$. Bounds that are sums of such terms are treated term by term; the two terms $(\beta,\gamma)=(\frac19,\frac12)$ and $(\frac16,\frac38)$ of the first bound in Lemma~\ref{lem:combine} give Proposition~\ref{prop:boundary}.
\end{remark}

\begin{corollary}\label{cor:boundary}
Let $0<H\le\frac14$. Then each of the four sums
\begin{alignat*}{2}
\mathcal B^+_1(R)&=\sum_{R^2<n\le(R+H)^2}\Np(n)f(\sqrt n),&\qquad \mathcal B^+_4(R)&=\sum_{R^2<4n\le(R+H)^2}\Np(4n)f(\sqrt{4n}),\\
\mathcal B^-_1(R)&=\sum_{R^2<n\le(R+H)^2}\Nm(n)f(\sqrt n),&\qquad \mathcal B^-_4(R)&=\sum_{R^2<4n\le(R+H)^2}\Nm(4n)f(\sqrt{4n})
\end{alignat*}
equals $c\int_R^{R+H}t^2f(t)\,dt$, with $c=\frac{\pi^2}6$, $\frac{\pi^2}{12}$, $\frac\pi6$ and $\frac\pi{12}$ respectively, up to an error bounded by the error term of Proposition~\ref{prop:boundary}. More precisely, any bound for $\mathcal B(R)-c(\sigma,\mathcal A)\int_R^{R+H}t^2f(t)\,dt$ that holds for all balanced families also bounds these errors. For $H=R^{-13/21}$ the errors are $O_\eps(R^{55/42+\eps})$.
\end{corollary}

\begin{proof}
By \eqref{eq:Lform}, writing $n=d^2a$ (resp.\ $4n=d^2a$), each sum is $c_\pm\mathcal B(R)$ for a suitable family, with $c_+=1$ and $c_-=1/\pi$, and then $c=c_\pm c(\sigma,\mathcal A)$:
\begin{itemize}
\item for $\mathcal B^+_1$: $\sigma=1$ and $\mathcal A_d=\{0,1,4,5\}$;
\item for $\mathcal B^+_4$: $\sigma=1$, $\mathcal A_d=\{0,4\}$ for odd $d$ and $\{0,1,4,5\}$ for even $d$;
\item for $\mathcal B^-_1$: $\sigma=-1$ and $\mathcal A_d=\{0,3,4,7\}$;
\item for $\mathcal B^-_4$: $\sigma=-1$, $\mathcal A_d=\{0,4\}$ for odd $d$ and $\{0,3,4,7\}$ for even $d$.
\end{itemize}
All four families are balanced. If $\mathcal A_d$ consists of all classes with $\sigma a\equiv0,1\pmod4$, then $\delta_d(b)=\frac12\varphi(b)/b$ for every $b$; if $\mathcal A_d=\{0,4\}$, then $\delta_d(b)=\frac14\varphi(b)/b$ for odd $b$ and $\delta_d(b)=0$ for even $b$. Since $\sum_b\varphi(b)b^{-3}=\zeta(2)/\zeta(3)$ and $\sum_{b\ \mathrm{odd}}\varphi(b)b^{-3}=\frac67\zeta(2)/\zeta(3)$, \eqref{eq:cconst} gives $c(\sigma,\mathcal A)=\zeta(2)$ for the first and third families, and
\[
c(\sigma,\mathcal A)=\frac{\zeta(2)}{\zeta(3)}\Big(\sum_{d\ \mathrm{even}}\frac1{d^3}+\frac37\sum_{d\ \mathrm{odd}}\frac1{d^3}\Big)=\zeta(2)\Big(\frac18+\frac38\Big)=\frac{\zeta(2)}2
\]
for the second and fourth. Multiplying by $c_+=1$ or $c_-=1/\pi$ gives, in the stated order, $c=\frac{\pi^2}6$, $\frac{\pi^2}{12}$, $\frac\pi6$ and $\frac\pi{12}$. With $H=R^{-13/21}$ the error terms of Proposition~\ref{prop:boundary} are $R^{1+\eps}$, $R^{82/63+\eps}$ and $R^{55/42+\eps}$, and $82/63<55/42$.
\end{proof}

\section{Exponential sums over class numbers}\label{sec:expsums}

In this section we will bound the dual sums in the summation formulas of Section~\ref{ssec:shintani}, which are exponential sums over class numbers. Let $a(n)$ be one of $\Np(n)$, $\Np(4n)$, $\Nm(n)$, $\Nm(4n)$, let $v$ be smooth and supported in $(\frac12,2)$, and for $\rho\ge2$ and $M\ge1$ put
\[
V_M(\rho;v)=\sum_{n\ge1}a(n)\,v(n/M)\,e(\rho\sqrt n).
\]

\begin{proposition}\label{prop:expsum}
Let $(\kappa,\lambda)$ be an exponent pair. Then for $\rho\ge2$ and $1\le M\le 100\rho^2$,
\[
V_M(\rho;v)\ll_\eps(\rho M)^\eps\big(M^{5/4}+\rho^{(1+\kappa)/4}M^{7/8+\lambda/4-\kappa/8}\big)
\]
uniformly, where the implied constant depends on $\eps$, the pair and finitely many derivatives of $v$.
\end{proposition}

The proof has two steps. First, the lattice point arguments of \cite{CI1}, \cite{CI2} and \cite{CU} reduce $V_M$ to a double sum (see Lemma~\ref{lem:reduction}). Second, Lemma~\ref{lem:BK41}, a version of \cite[Lemma 4.1]{BK} that is uniform in the summation intervals, estimates that double sum.

\begin{lemma}\label{lem:reduction}
Let $0<\delta\le\frac1{10}$, let $a(n)$ be as above, let $I\subseteq[M,2M]$ be an interval and $\rho\in\mathbb R$. Then
\[
\Big|\sum_{n\in I}a(n)e(\rho\sqrt n)\Big|^2\ll_{\eps}M^{5/2+3\delta+\eps}+M^{1+3\delta+\eps}\sum_{1\le w\le W}\Big|\sum_{u\in I_w}e\big(\rho(\sqrt u-\sqrt{u+w})\big)\Big|,
\]
where $W=16M^{1-\delta}$ and $I_w=I\cap(I-w)$. For $a(n)=\Np(4n)$ or $\Nm(4n)$ the same holds with $4M$, $4I$ and $\rho/2$ in place of $M$, $I$ and $\rho$.
\end{lemma}

\begin{proof}
First consider $a(n)=\Np(n)$; we follow \cite[Proof of Proposition 3.3]{CU}. Put $E(m)=e(\rho\sqrt m)\mathbf 1_{m\in I}$. Since $\Np(n)=0$ for square $n$, Lemma~\ref{lem:siegel} gives
\[
\sum_{n\in I}\Np(n)e(\rho\sqrt n)=\sum_{\substack{a>0,\ a+c<|b|/2\\ b^2-4ac\ne\sq}}\mu(a,b,c)E(b^2-4ac),
\]
and $c\ne0$ in every term. The length $\mu(a,b,c)$ is defined, with $\mu\ll M^\eps$, also when $b^2-4ac=\sq$, and the triples with $c\ne0$, $a+c<|b|/2$ and $b^2-4ac=\sq\in I$ number $O(M^{1+\eps})$ \cite[Proof of Proposition 3.3]{CU}. Adding them, and summing by parts in $b$ over the $O(1)$ intervals on which $\mu(a,\cdot,c)$ is monotone,
\[
\Big|\sum_{n\in I}\Np(n)e(\rho\sqrt n)\Big|\ll M^\eps\sum_{1\le|4ac|\le2M}\ \max_{I'}\Big|\sum_{b\in I'}E(b^2-4ac)\Big|+M^{1+\eps},
\]
with $I'\subseteq[-2\sqrt M,2\sqrt M]$. Each $n=4ac\ne0$ arises $O(M^\eps)$ times, so by \cite[Lemma 7.3]{GK} there is $\theta\in\mathbb R$ with
\begin{equation}\label{eq:red1}
\Big|\sum_{n\in I}\Np(n)e(\rho\sqrt n)\Big|\ll M^{2\eps}\sum_{1\le|n|\le2M}\Big|\sum_{|b|\le2\sqrt M}e(\theta b)E(b^2-n)\Big|+M^{1+\eps}.
\end{equation}
Split $|b|\le2\sqrt M$ into $O(M^\delta)$ intervals $\mathcal J$ of length $\le M^{1/2-\delta}$. By Cauchy's inequality in $\mathcal J$ and $n$, extending the $n$-sum to $\Z$, the square of the sum over $n$ in \eqref{eq:red1} is
\begin{equation}\label{eq:red2}
\ll M^{1+\delta}\sum_{\mathcal J}\ \sum_{b_1,b_2\in\mathcal J}\Big|\sum_{n\in\Z}E(b_1^2-n)\overline{E(b_2^2-n)}\Big|.
\end{equation}
The terms $|b_1|=|b_2|$ contribute $O(M^{3/2})$ to the double sum. For $|b_1|<|b_2|$ put $u=b_1^2-n$ and $w=b_2^2-b_1^2$; then $1\le w\le|b_2-b_1|(|b_1|+|b_2|)\le4M^{1-\delta}$ and
\[
\sum_{n\in\Z}E(b_1^2-n)\overline{E(b_2^2-n)}=\sum_{u\in I\cap(I-w)}e\big(\rho(\sqrt u-\sqrt{u+w})\big).
\]
Each $w$ arises from $O(M^\eps)$ pairs $(b_1,b_2)$, and the terms $|b_1|>|b_2|$ are the conjugates. With \eqref{eq:red1} and \eqref{eq:red2},
\[
\Big|\sum_{n\in I}\Np(n)e(\rho\sqrt n)\Big|^2\ll M^{5/2+\delta+5\eps}+M^{1+\delta+5\eps}\sum_{1\le w\le4M^{1-\delta}}\Big|\sum_{u\in I_w}e\big(\rho(\sqrt u-\sqrt{u+w})\big)\Big|.
\]

For $a(n)=\Nm(n)$ use Lemma~\ref{lem:gauss}. The corrections at $n=3k^2,4k^2$ contribute $O(M^{1/2})$, and the forms with $|b|=a$ or $a=c$ number $O(M^{1+\eps})$. For the others there is no weight, $3a^2\le4ac-b^2\le2M$, so $a\le\sqrt M$, and for fixed $(a,c)$ the variable $b$ runs over an interval. The argument is the same, with $n=4ac\in[1,3M]$, $E(n-b^2)$ in place of $E(b^2-n)$, and $u=n-b_2^2$; the correlation sums become
\[
\sum_{u\in I\cap(I-w)}e\big(\rho(\sqrt{u+w}-\sqrt u)\big),
\]
the conjugates of those above.

For $a(n)=\Np(4n)$ or $\Nm(4n)$, the sum is over $n'=4n\in4I$, with phase $e(\frac\rho2\sqrt{n'})$. In the lattice representations $n'\equiv0\pmod4$ means $2\mid b$, and $\mathbf 1_{2\mid b}=\frac12(1+e(b/2))$; the phase $e(b/2)$ is absorbed into $e(\theta b)$ (see \cite[Proof of Proposition 3.3]{CU}).
\end{proof}

\begin{lemma}\label{lem:BK41}
Let $(\kappa,\lambda)$ be an exponent pair, $0<\delta\le\frac1{10}$, $\rho\ge2$ and $1\le M\le10^4\rho^2$. Let $1\le W\le16M^{1-\delta}$, and for $1\le w\le W$ let $I_w\subseteq[M,2M]$ be arbitrary intervals. Then
\[
\sum_{1\le w\le W}\Big|\sum_{u\in I_w}e\big(\rho(\sqrt u-\sqrt{u+w})\big)\Big|\ll_{\eps}(\rho M)^\eps\big(M^{3/2}+\rho^{(1+\kappa)/2}M^{3/4+\lambda/2-\kappa/4}\big),
\]
the implied constant depending on $\eps$, $\delta$ and the exponent pair.
\end{lemma}

\begin{proof}
We follow \cite[Proof of Lemma 4.1]{BK}, which starts from \cite[(3.4)]{CI1}, and add a completion step for the intervals $I_w$. Let $\varphi(u,w)=\rho(\sqrt u-\sqrt{u+w})$ and, for dyadic $D\le16M^{1-\delta}$,
\[
S_D=\sum_{D\le w<2D}\Big|\sum_{u\in I_w}e(\varphi(u,w))\Big|.
\]
It suffices to show that
\begin{equation}\label{eq:SD}
S_D^2\ll(\rho M)^\eps\big(M^3+\rho^{1+\kappa}M^{3/2+\lambda-\kappa/2}\big).
\end{equation}
We treat small derivatives directly; otherwise we complete the intervals $I_w$, apply Poisson summation with stationary phase and Cauchy's inequality in $w$, and conclude with the exponent pair.
For $u\in[M,2M]$, $w\in[D,2D)$ and $j\ge1$,
\begin{equation}\label{eq:phi-der}
\varphi_u=\frac\rho2\big(u^{-1/2}-(u+w)^{-1/2}\big)\asymp U:=\rho DM^{-3/2},\qquad|\partial_u^j\varphi|\asymp_jUM^{1-j},
\end{equation}
and $\varphi_u$ is decreasing in $u$ and increasing in $w$. Let $c_0>0$ be such that $U\le c_0$ implies $0<\varphi_u\le\frac12$. If $U\le c_0$, the Kusmin--Landau inequality \cite[Theorem 2.1]{GK} gives, for any intervals $I_w$,
\[
S_D\ll\frac DU=\frac{M^{3/2}}\rho\ll M^{3/2}.
\]
It remains to treat $U>c_0$. For integers $u\in[M,2M]$,
\[
\mathbf 1_{I_w}(u)=\int_0^1c_w(\theta)e(\theta u)\,d\theta,\qquad c_w(\theta)=\sum_{v\in I_w}e(-\theta v),
\]
where $|c_w(\theta)|\le K(\theta):=\min\{M+1,\|\theta\|^{-1}\}$
and $\int_0^1K\ll\log2M$. Hence
\begin{equation}\label{eq:completion}
S_D\ll\log(2M)\sup_{0\le\theta<1}S_D(\theta),\qquad S_D(\theta):=\sum_{D\le w<2D}\Big|\sum_{M\le u\le2M}e\big(\varphi(u,w)+\theta u\big)\Big|.
\end{equation}
Fix $\theta$ and $w$, and let $\psi(u)=\varphi(u,w)+\theta u$. Then $\psi'$ is decreasing, $\psi''\asymp-U/M$ and $|\psi^{(j)}|\asymp UM^{1-j}$ for $j=3,4$, uniformly in $\theta$, and \cite[Lemma 3.6]{GK} gives
\begin{equation}\label{eq:poisson}
\sum_{M\le u\le2M}e(\psi(u))=e(-\tfrac18)\sum_{s\in\mathcal S_w}A(s,w)\,e\big(g(s,w)\big)+O\big(\log(2+U)+(M/U)^{1/2}\big),
\end{equation}
where $\mathcal S_w=(\Z-\theta)\cap[\varphi_u(2M,w),\varphi_u(M,w)]$, $\alpha=\alpha(s,w)\in[M,2M]$ solves $\varphi_u(\alpha,w)=s$, and
\[
A(s,w)=|\varphi_{uu}(\alpha,w)|^{-1/2}\asymp(M/U)^{1/2},\qquad g(s,w)=\varphi(\alpha,w)-s\alpha.
\]
Indeed, for $k\in\Z$ the stationary point of $\psi(u)-ku$ solves $\varphi_u=k-\theta=s$, and the phase there is $g(s,w)$. These are the functions of \cite[(4.6)--(4.7)]{BK} with $y=w$, at the points $s\in\Z-\theta$; since \cite[(4.6)--(4.11)]{BK} treat $s\asymp U$ as a real variable, they apply unchanged. As $D\le16M$, $M\le10^4\rho^2$ and $U\le16\rho M^{-1/2}$, the error terms in \eqref{eq:poisson} contribute to $S_D(\theta)^2$
\[
\ll D^2\log^2(2+U)+\frac{DM^{5/2}}\rho\ll_\eps(\rho M)^\eps M^3.
\]
By Cauchy's inequality in $w$, the square of the main part of $S_D(\theta)$ is at most
\[
D\sum_{D\le w<2D}\ \sum_{s,\,s+z\in\mathcal S_w}A(s,w)A(s+z,w)\,e\big(G(s,w,z)\big),\qquad G(s,w,z)=g(s+z,w)-g(s,w),
\]
with $s$ in an interval of length $O(U)$ and $z\in\Z$, $|z|\ll U$. Since $\varphi_u$ is increasing in $w$, the set $\mathcal W(s,z)$ of $w\in[D,2D)$ with $s,s+z\in\mathcal S_w$ is an interval, and the last expression is
\begin{equation}\label{eq:cauchy}
\le D\sum_s\sum_z\Big|\sum_{w\in\mathcal W(s,z)}A(s,w)A(s+z,w)\,e\big(G(s,w,z)\big)\Big|.
\end{equation}
By the scaling $\alpha=a_0\mathcal A(\eta)$ of \cite[(4.8)--(4.9)]{BK}, $\partial_w^jA(s,w)\ll_j(M/U)^{1/2}D^{-j}$ for $j\le1$, uniformly in $s\asymp U$; by partial summation the inner sum in \eqref{eq:cauchy} is
\[
\ll\frac MU\max_{\mathcal W'\subseteq\mathcal W(s,z)}\Big|\sum_{w\in\mathcal W'}e\big(G(s,w,z)\big)\Big|.
\]
The terms $z=0$ contribute $D\cdot O(U)\cdot O(M/U)\cdot O(D)=O(M^3)$. Since $G(s,w,z)=-G(s+z,w,-z)$ and $\mathcal W(s,z)=\mathcal W(s+z,-z)$, the terms $z<0$ are the conjugates of those with $z>0$. For $z\asymp T$, $1\le T\ll U$, \cite[(4.7)--(4.11)]{BK} give
\[
-G(s,w,z)=\Lambda_{s,z}w^{2/3}+u_{s,z}(w),\qquad\Lambda_{s,z}=\int_s^{s+z}\Big(\frac\rho{4v}\Big)^{2/3}dv>0,\qquad\Lambda_{s,z}D^{2/3}\asymp TM,
\]
using $\eta\asymp D/M\ll M^{-\delta}$, with $u_{s,z}$ satisfying \cite[(4.10)]{BK}. Hence, for $M\gg_{\delta,\kappa,\lambda}1$, the stability hypothesis \cite[(4.2)]{BK} holds, and \cite[(4.12)]{BK} gives, for every subinterval $\mathcal W'$,
\[
\sum_{w\in\mathcal W'}e\big(G(s,w,z)\big)\ll(TM)^\kappa D^{\lambda-\kappa}+\frac D{TM}.
\]
There are $O(U)$ values of $s$ and $O(T)$ values of $z$, and, as in \cite[p.~15]{BK},
\[
UT(TM)^\kappa D^{\lambda-\kappa}\ll\rho^{2+\kappa}M^{-1+\lambda-\kappa/2},\qquad UT\cdot\frac D{TM}\ll\rho M^{-1/2}.
\]
Multiplying by $D\cdot M/U=\rho^{-1}M^{5/2}$ and summing over dyadic $T$, the terms $z\ne0$ contribute
\[
\ll(\rho M)^\eps\big(\rho^{1+\kappa}M^{3/2+\lambda-\kappa/2}+M^2\big).
\]
Hence $S_D(\theta)$ satisfies \eqref{eq:SD} uniformly in $\theta$, and \eqref{eq:SD} follows by \eqref{eq:completion}. In \cite{BK} the hypothesis $M\le\rho^2$ is used only for the error term in \eqref{eq:poisson}; here $M\le10^4\rho^2$ suffices.
\end{proof}

\begin{proof}[Proof of Proposition~\ref{prop:expsum}]
For $2\le\rho<4$ we have $M\le1600$, and the bound is trivial; let $\rho\ge4$. By partial summation,
\[
|V_M(\rho;v)|\ll\max_{I\subseteq[M/2,2M]}\Big|\sum_{n\in I}a(n)e(\rho\sqrt n)\Big|,
\]
and each such $I$ is the union of at most two intervals as in Lemma~\ref{lem:reduction}, with $M/2$ or $M$ in place of $M$. Apply Lemma~\ref{lem:reduction} and then Lemma~\ref{lem:BK41}, with $M/2$, $M$, $2M$ or $4M$ in place of $M$ and $\rho$ or $\rho/2\ge2$ in place of $\rho$; the condition $M\le10^4\rho^2$ holds. Finally take $\delta$ small in terms of $\eps$.
\end{proof}

\begin{remark}\label{rem:pairs}
With the exponent pair $(\kappa,\lambda)=(\frac12,\frac12)$, Proposition~\ref{prop:expsum} gives $V_M\ll(\rho M)^\eps(M^{5/4}+\rho^{3/8}M^{15/16})$. The term $\rho^{3/8}M^{15/16}$ occurs in the bounds of \cite[Lemma 3.1]{CI1}, \cite[Lemma 4.1]{CI2} and \cite[Proposition 3.3]{CU}; the last also contains $\rho^{1/8}M^{17/16}$, which for $M\le100\rho^2$ is $O(\rho^{3/8}M^{15/16})$. With Bourgain's exponent pair $(\frac{13}{84}+\eta,\frac{55}{84}+\eta)$ \cite[Theorem 6]{Bou} it gives $V_M\ll(\rho M)^\eps(M^{5/4}+\rho^{97/336}M^{685/672})$, as in \cite[(4.13)]{BK}; here the extra factor $\rho^{\eta/4}M^{\eta/8}$ is absorbed into $(\rho M)^\eps$ by taking $\eta$ small.
\end{remark}

\begin{corollary}\label{cor:Vbound}
Let $R\ge16$, $\frac R2\le\rho\le2R$ and $1\le M\le R^2$. Then $V_M(\rho;v)/M\ll_\eps R^\eps\big(M^{1/4}+R^{17/55}\big)$.
\end{corollary}

\begin{proof}
By Remark~\ref{rem:pairs},
\begin{align*}
\frac{V_M}M&\ll R^\eps\Big(M^{1/4}+\min\big\{R^{3/8}M^{-1/16},R^{97/336}M^{13/672}\big\}\Big)\\
&\le R^\eps\Big(M^{1/4}+\big(R^{3/8}M^{-1/16}\big)^{13/55}\big(R^{97/336}M^{13/672}\big)^{42/55}\Big),
\end{align*}
and $\frac38\cdot\frac{13}{55}+\frac{97}{336}\cdot\frac{42}{55}=\frac{17}{55}$, $-\frac1{16}\cdot\frac{13}{55}+\frac{13}{672}\cdot\frac{42}{55}=0$.
\end{proof}

\section{Smoothed sums}\label{sec:smoothed}

We now combine the summation formulas of Section~\ref{ssec:shintani} with the bounds of Section~\ref{sec:expsums} to evaluate the smoothed sums. Since our weight is smooth, Shintani's summation formula applies directly, without the mollification used in \cite{CU}. Together with the exponent pair $(\frac12,\frac12)$ in Remark~\ref{rem:pairs}, this means that the term $N^{11/16+\eps}\Delta^{1/8}$ of \cite[Proposition 4.1]{CU} does not arise; with $\Delta=H=R^{-13/21}$ it would in any case be $R^{109/84+\eps}$, below $R^{55/42}$. We begin with the dual sums.

\begin{lemma}\label{lem:dual}
Let $R^{-1}\le H\le\frac14$. Let $a(n)$ be one of $\Np(4n)$, $\Np(n)$, $\Nm(4n)$, $\Nm(n)$, let $T\in\{I,J\}$ be as in Lemma~\ref{lem:transforms}, and let $\lambda_0\in\{1,\frac12\}$. Suppose that, for some $p_1\ge0$ and $0\le p_2<\frac12$,
\begin{equation}\label{eq:Vhyp}
V_M(\rho;v)/M\ll_\eps R^\eps\big(M^{1/4}+R^{p_1}M^{p_2}\big)
\end{equation}
for $\frac R2\le\rho\le2R$ and $1\le M\le R^2$, for every smooth $v$ supported in $(\frac12,2)$, the implied constant depending on $\eps$ and $v$. Then
\[
\sum_{n\ge1}\frac{a(n)}{\sqrt n}\,T(\lambda_0\sqrt n)\ll_\eps R^\eps\big(RH^{-1/2}+R^{1+p_1}H^{-2p_2}+R^{-1}H^{-3}\big).
\]
By Corollary~\ref{cor:Vbound}, \eqref{eq:Vhyp} holds with $(p_1,p_2)=(\frac{17}{55},0)$; for $H=R^{-13/21}$ the sum is then $O_\eps(R^{55/42+\eps})$.
\end{lemma}

\begin{proof}
Let $w\in C_0^\infty((\frac12,2))$ with $\sum_{j\in\Z}w(2^{-j}x)=1$ for $x>0$, let $w_\alpha(x)=x^{-\alpha}w(x)$, and write the sum as $\sum_{M=2^j\ge1}E_M$, where
\[
E_M=\sum_n\frac{a(n)}{\sqrt n}w\Big(\frac nM\Big)T(\lambda_0\sqrt n).
\]
The terms $O_A(\xi^{-A})$ of Lemma~\ref{lem:transforms} contribute $O(1)$ in total. Inserting the representations of Lemma~\ref{lem:transforms} and interchanging $\int_R^{R+H}$ with $\sum_n$, we see that $E_M$ is a linear combination of the integrals
\[
\int\frac{V_M(\pm\lambda_0t;w_1)}M\,tf'(t)\,dt,\qquad\int\frac{V_M(\pm\lambda_0t;w_{3/2})}{M^{3/2}}\,f'(t)\,dt,\qquad\int\frac{V_M(\pm\lambda_0t;w_{3/2})}{M^{3/2}}\,F''(t)\,dt,
\]
the first two from the first representation and the third from the second. Since $\|tf'\|_1\ll R$, $\|f'\|_1\ll1$ and $\|F''\|_1\ll R/H$, as in \cite[(4.18)]{BK},
\begin{equation}\label{eq:EM}
|E_M|\ll\frac RM\min\big\{1,(H\sqrt M)^{-1}\big\}\sup_{\lambda_0R\le\rho\le\lambda_0(R+H)}\ \max_{\alpha\in\{1,3/2\}}|V_M(\rho;w_\alpha)|+O(1).
\end{equation}
For $M\le R^2$ we have $\frac R2\le\rho\le2R$, and \eqref{eq:Vhyp} and \eqref{eq:EM} give, splitting at $M=H^{-2}\le R^2$ as in \cite[(4.19)]{BK},
\begin{align*}
\sum_{M\le R^2}|E_M|&\ll R^{1+\eps}\Big(\sum_{M\le H^{-2}}\big(M^{1/4}+R^{p_1}M^{p_2}\big)+H^{-1}\sum_{M>H^{-2}}\big(M^{-1/4}+R^{p_1}M^{p_2-1/2}\big)\Big)\\
&\ll R^{1+\eps}\big(H^{-1/2}+R^{p_1}H^{-2p_2}\log R\big),
\end{align*}
since $p_2<\frac12$. For $M>R^2$, $|T(\xi)|\ll RH^{-3}\xi^{-4}$ and $a(n)\ll n^{1/2+\eps}$ (by \eqref{eq:Lform}) give
\[
|E_M|\ll RH^{-3}M^{-1+\eps},\qquad\sum_{M>R^2}|E_M|\ll R^{-1+\eps}H^{-3}.
\]
For $(p_1,p_2)=(\frac{17}{55},0)$ and $H=R^{-13/21}$ the three terms are $R^{55/42}$, $R^{72/55}$ and $R^{6/7}$.
\end{proof}

\begin{proposition}\label{prop:smoothed}
Let $R^{-1}\le H\le\frac14$, and let $\mathcal D\ge R$ bound the absolute values of the sums in Lemma~\ref{lem:dual}, for all choices of $a(n)$, $T$ and $\lambda_0$ there. Then
\begin{align*}
\sum_{n\ge1}\Np(n)f(\sqrt n)&=M_1(g)+O_\eps(R^\eps\mathcal D),\\
\sum_{n\ge1}\Np(4n)f(\sqrt{4n})&=M_4(g)+O_\eps(R^\eps\mathcal D),\\
\sum_{n\ge1}\Nm(n)f(\sqrt n)&=\frac\pi6\int_0^\infty tg(t)\,dt-\frac12\int_0^\infty g(t)\,dt+O_\eps(R^\eps\mathcal D),\\
\sum_{n\ge1}\Nm(4n)f(\sqrt{4n})&=\frac\pi{12}\int_0^\infty tg(t)\,dt-\frac14\int_0^\infty g(t)\,dt+O_\eps(R^\eps\mathcal D).
\end{align*}
By Lemma~\ref{lem:dual}, for $H=R^{-13/21}$ one may take $\mathcal D\ll_\eps R^{55/42+\eps}$.
\end{proposition}

\begin{proof}
Since $g(\sqrt n)/\sqrt n=f(\sqrt n)$, the left sides are those of Lemmas~\ref{lem:sum-pos} and~\ref{lem:sum-neg}. As $\tilde g(\sqrt{4n})=J(\sqrt n)$, $\tilde g(\sqrt n)=J(\sqrt n/2)$, $\breve g(\sqrt{4n})=I(\sqrt n)$ and $\breve g(\sqrt n)=I(\sqrt n/2)$, the dual sums there are
\[
\sum_n\frac{\Np(4n)}{\sqrt n}J(\sqrt n),\quad \sum_n\frac{\Np(n)}{\sqrt n}J\Big(\frac{\sqrt n}2\Big),\quad -\sum_n\frac{\Nm(4n)}{\sqrt n}I(\sqrt n),\quad -\sum_n\frac{\Nm(n)}{\sqrt n}I\Big(\frac{\sqrt n}2\Big),
\]
all $O(\mathcal D)$. As $0\le g(t)\le t$, the integrals against $(e^{2\pi t}-1)^{-1}$ and $(e^{\pi t}-1)^{-1}$ are $O(1)$. It remains to bound $E_1(g)$ and $E_4(g)$. By Lemma~\ref{lem:transforms} with $\xi=x/2$,
\begin{equation}\label{eq:gtilde}
\tilde g(x)\ll x^{-2}+R\min\{x^{-1},H^{-1}x^{-2}\}\qquad(x\ge1).
\end{equation}
By \eqref{eq:gtilde} and $\alpha_n,\beta_n\ll\log2n$, splitting at $n=H^{-1}$,
\[
\sum_n\log(2n)|\tilde g(n)|\ll1+R\log^2(2/H).
\]
Since $g=0$ on $[0,\frac14]$ and $g(t)\le t$, for $c>0$
\[
\int_0^\infty g(t)e^{-c\sqrt n\,t}\,dt\le\int_{1/4}^\infty te^{-c\sqrt n\,t}\,dt\ll_ce^{-c\sqrt n/4}\big(n^{-1/2}+n^{-1}\big),
\]
so, by $\Nm(n)\ll n^{1/2+\eps}$, the sums over $\Nm(4n)$ and $\Nm(n)$ in $E_1,E_4$ are $O(1)$. Here the vanishing of $g$ on $[0,\frac14]$ is used; the bound $g(t)\le t$ alone would give a divergent sum. For $c\in\{1,2\}$ and $0<t\le1$ we have $cn/t\ge1$, and by \eqref{eq:gtilde}, splitting at $n=t/(cH)$,
\[
\sum_n|\tilde g(cn/t)|\ll t^2+Rt\big(1+\log^+(t/H)\big),\qquad\int_0^1\sum_n|\tilde g(cn/t)|\,\frac{dt}{t(1+t)}\ll1+R\log(2/H),
\]
where $\log^+x=\max\{0,\log x\}$. Hence $E_1(g),E_4(g)\ll R\log^2R\ll R^\eps\mathcal D$, as $\mathcal D\ge R$ and $H\ge R^{-1}$.
\end{proof}

\section{Transfer principle and proofs of the theorems}\label{sec:proofs}

Subtracting the boundary sums of Corollary~\ref{cor:boundary} from the smoothed sums of Proposition~\ref{prop:smoothed} gives asymptotic formulas for the sums of $N_\pm$. The theorems then follow by M\"{o}bius inversion. We state this in a form that applies to each of our three sets of hypotheses.

\subsection{Transfer to the sharp sums}

\begin{proposition}\label{prop:Nsums}
Let $\frac12<\vartheta<1$. Suppose that for every sufficiently large $R$ there is an $H\in[R^{-1},\frac14]$ such that the four boundary sums of Corollary~\ref{cor:boundary} equal their main terms up to $O_\eps(R^{2\vartheta+\eps})$, and the sums in Lemma~\ref{lem:dual} are $O_\eps(R^{2\vartheta+\eps})$ for all choices of $a(n)$, $T$ and $\lambda_0$. Then, for $x\ge1$,
\begin{align*}
\sum_{n\le x}\Np(n)&=\frac{\pi^2}{18}x^{3/2}-\frac x2\log x+\frac x2\Big(1-\frac{\zeta'(2)}{\zeta(2)}-\log2\pi\Big)+O_\eps(x^{\vartheta+\eps}),\\
\sum_{n\le x}\Np(4n)&=\frac{2\pi^2}{9}x^{3/2}-x\log x+x\Big(1-\frac{\zeta'(2)}{\zeta(2)}-\log2\pi-\frac53\log2\Big)+O_\eps(x^{\vartheta+\eps}),\\
\sum_{n\le x}\Nm(n)&=\frac{\pi}{18}x^{3/2}-\frac x4+O_\eps(x^{\vartheta+\eps}),\\
\sum_{n\le x}\Nm(4n)&=\frac{2\pi}{9}x^{3/2}-\frac x2+O_\eps(x^{\vartheta+\eps}).
\end{align*}
\end{proposition}

\begin{proof}
Bounded $x$ are trivial. Let $R=\sqrt x$ for the first and third sums and $R=2\sqrt x$ for the second and fourth, so that $\sqrt{4n}\le R$ if and only if $n\le x$, and let $H$ be as in the hypothesis. Since $f(\sqrt n)=1$ for $n\le R^2$,
\[
\sum_{n\le x}\Nm(n)=\sum_{n\ge1}\Nm(n)f(\sqrt n)-\mathcal B^-_1(R),
\]
and similarly for the other three sums. By Proposition~\ref{prop:smoothed} with $\mathcal D\ll_\eps R^{2\vartheta+\eps}$ (admissible, as $2\vartheta>1$) and the hypothesis on Corollary~\ref{cor:boundary}, the error terms are $O(R^{2\vartheta+\eps})=O(x^{\vartheta+\eps})$. Since $g(t)=t$ on $[\frac12,R]$,
\[
c\int_0^\infty tg(t)\,dt-c\int_R^{R+H}t^2f(t)\,dt=c\int_0^Rt^2\,dt+O(1)=\frac c3R^3+O(1),
\]
which gives $\frac{\pi^2}{18}x^{3/2}$, $\frac{\pi^2}{36}(2\sqrt x)^3=\frac{2\pi^2}9x^{3/2}$, $\frac\pi{18}x^{3/2}$ and $\frac{2\pi}9x^{3/2}$. For $\Nm$, $\int_0^\infty g=\frac{R^2}2+O(R)$ gives $-\frac x4$ and $-\frac14\cdot\frac{4x}2=-\frac x2$. The other terms of $M_1(g)$ and $M_4(g)$ depend only on $g|_{[0,R+H]}$ and $g|_{\Z}$; replacing $g$ by $t\,\mathbf 1_{[0,R]}(t)$ changes them by $O(R^{1+\eps})$, and we are in the situation of \cite[Proof of Proposition 4.1]{CU}. With
\[
\sum_{n\le y}n\log n=\frac{y^2}2\log y-\frac{y^2}4+O(y\log y),\qquad\sum_d\frac{\Lambda(d)}{d^2}=-\frac{\zeta'(2)}{\zeta(2)},
\]
one finds for $R^2=x$, up to $O(R^{1+\eps})$,
\begin{align*}
\int_0^\infty g(t)\log(2\pi t)\,dt&=\frac x4\log x+\frac{2\log2\pi-1}4x,\\
\sum_ng(n)\log n-\sum_d\frac{\Lambda(d)}d\sum_ng(dn)&=\frac x4\log x-\frac x4+\frac{\zeta'(2)}{2\zeta(2)}x.
\end{align*}
This gives the first formula. For $R^2=4x$ the corresponding terms of $M_4(g)$ are
\[
-\frac x2\log x-x\log4\pi+\frac x2,\qquad -\frac x2\log x+\frac x2-\frac{\zeta'(2)}{\zeta(2)}x,\qquad -\frac23x\log2,
\]
the last coming from $-\log2\sum_k2^{-k}\sum_ng(2^kn)=-\log2\cdot\frac{R^2}2\sum_{k\ge1}4^{-k}$. Their sum gives the second formula.
\end{proof}

\begin{corollary}\label{cor:transfer}
Under the hypothesis of Proposition~\ref{prop:Nsums}, the four formulas of Theorems~\ref{thm:neg} and~\ref{thm:pos} hold with error terms $O_\eps(N^{\vartheta+\eps})$ in place of $O_\eps(N^{55/84+\eps})$.
\end{corollary}

\begin{proof}
By Lemma~\ref{lem:mobius},
\[
\sum_{n\le N}h(-n)=\sum_{k\le\sqrt N}\mu(k)\sum_{m\le N/k^2}\Nm(m)+O(1),
\]
the $O(1)$ accounting for $\delta_{-3}$ and $\delta_{-4}$. Insert Proposition~\ref{prop:Nsums} and complete the sums over $k$, at a cost $O(N^{1/2})$. Since $\vartheta>\frac12$,
\[
\sum_k\Big(\frac N{k^2}\Big)^{\vartheta+\eps}\ll N^{\vartheta+\eps},\qquad\sum_k\frac{\mu(k)}{k^3}=\frac1{\zeta(3)},\qquad\sum_k\frac{\mu(k)}{k^2}=\frac6{\pi^2},
\]
and the main term is $\frac{\pi}{18\zeta(3)}N^{3/2}-\frac{3}{2\pi^2}N$. Similarly, with $A_j(y)=\sum_{m\le y}\Nm(jm)$,
\[
\sum_{n\le N}h(-4n)=\sum_{k\ \mathrm{odd}}\mu(k)\,(A_4-A_1)\Big(\frac N{k^2}\Big)+O(1),\qquad(A_4-A_1)(y)=\frac\pi6y^{3/2}-\frac y4+O(y^{\vartheta+\eps}),
\]
and $\sum_{k\ \mathrm{odd}}\mu(k)k^{-3}=\frac8{7\zeta(3)}$, $\sum_{k\ \mathrm{odd}}\mu(k)k^{-2}=\frac8{\pi^2}$ give $\frac{4\pi}{21\zeta(3)}N^{3/2}-\frac2{\pi^2}N$. For positive discriminants the argument is the same. The main terms are evaluated in \cite[Section 4]{CU}, using $\sum_k\mu(k)k^{-2}\log k=\zeta'(2)/\zeta(2)^2$ and its analogue over odd $k$, and are those of \cite[Theorems 1.1, 1.2]{CU}; for instance, with $K=1-\zeta'(2)/\zeta(2)-\log2\pi$,
\begin{multline*}
\sum_{k}\mu(k)\Big(-\frac{N}{2k^2}\log\frac N{k^2}+\frac N{2k^2}K\Big)
=-\frac3{\pi^2}N\log N+\frac{\zeta'(2)}{\zeta(2)^2}N+\frac3{\pi^2}NK
=-\frac3{\pi^2}(C+\log N)N.
\end{multline*}
\end{proof}

\subsection{Unconditional estimates}

\begin{proof}[Proof of Theorems~\ref{thm:neg} and~\ref{thm:pos}]
By Corollary~\ref{cor:boundary} and Lemma~\ref{lem:dual} with $H=R^{-13/21}$, the hypothesis of Proposition~\ref{prop:Nsums} holds with $\vartheta=\frac{55}{84}$. Apply Corollary~\ref{cor:transfer}.
\end{proof}

The proofs of Theorem~\ref{thm:cond} and of Remark~\ref{rem:grh} below verify the hypothesis of Proposition~\ref{prop:Nsums} with other character-sum and exponential-sum bounds and other choices of H. Table~\ref{tab:opt} summarizes the three cases.

\begin{table}[ht]
\centering
\footnotesize
\begin{tabular}{@{}>{\raggedright\arraybackslash}p{2.4cm}>{\raggedright\arraybackslash}p{3.3cm}lllc@{}}
\toprule
Hypothesis & Boundary estimate & $(p_1,p_2)$ & $H$ & Dominant terms & $\vartheta$\\
\midrule
None & Lemma~\ref{lem:combine}, first bound & $(\frac{17}{55},0)$ & $R^{-13/21}$ & $R^{7/6}(RH)^{3/8}$, $RH^{-1/2}$ & $\frac{55}{84}$\\[4pt]
No zeros in $\mathrm{Re}\,s>\frac78$ & Lemma~\ref{lem:combine}, second bound; Lemma~\ref{lem:Lapprox} & $(p,0)$ & $R^{-3/5}$ & $R^{1+p}$ & $\frac{1+p}2$\\[4pt]
Lindelöf (GRH) & $(\beta,\gamma)=(0,\frac12)$ & $(\frac{97}{336},\frac{13}{672})$ & $R^{-71/181}$ & $R^{\frac32}H^{\frac12}$, $R^{\frac{433}{336}}H^{-\frac{13}{336}}$ & $\frac{118}{181}$\\
\bottomrule
\end{tabular}
\vspace{6pt}
\caption{The three applications of Proposition~\ref{prop:Nsums}. Here $(\beta,\gamma)$ refers to Remark~\ref{rem:general-char}, $(p_1,p_2)$ to \eqref{eq:Vhyp}, and $p=\frac{530939}{1743872}$, so that $\frac{1+p}2=\frac{2274811}{3487744}$. The dominant error terms are those equal to $R^{2\vartheta}$; all other terms are smaller.}\label{tab:opt}
\end{table}

\subsection{The zero-free hypothesis}

In this section we will give a proof of Theorem~\ref{thm:cond} by using the following approximation to $L(1,\chi)$. This approximation is similar to \cite[(6.1)]{BK}.

\begin{lemma}\label{lem:Lapprox}
Assume that no Dirichlet $L$-function has a zero in the half-plane $\mathrm{Re}\,s>\frac78$, and let $U\in C^\infty([0,\infty))$ be a smooth cut-off with $U=1$ on $[0,1]$ and $U=0$ on $[2,\infty)$. Then for every non-principal quadratic Dirichlet character $\chi$ modulo $k$ and every $Y\ge1$,
\[
L(1,\chi)=\sum_{m\ge1}\frac{\chi(m)}mU\Big(\frac mY\Big)+O_\eps\big((kY)^\eps Y^{-1/8}\big),
\]
the implied constant depending only on $\eps$ and $U$.
\end{lemma}

\begin{proof}
Let $\chi$ be induced by the primitive character $\chi^*$ of conductor $q\mid k$, let $0<\delta\le\frac1{16}$, and fix $t$; put $\mathcal L=\log(q(|t|+2))$. By hypothesis $L(s,\chi^*)$ has no zeros in the simply connected half-plane $\mathrm{Re}\,s>\frac78$, so it has there an analytic logarithm $\log L(s,\chi^*)$, normalized to agree with $\sum_p\sum_{j\ge1}\chi^*(p)^jj^{-1}p^{-js}$ for $\mathrm{Re}\,s>1$; and
\begin{equation}\label{eq:logL}
\begin{aligned}
|\log L(s,\chi^*)|&\le\log\zeta(\tfrac32)&&(\mathrm{Re}\,s\ge\tfrac32),\\
\mathrm{Re}\log L(s,\chi^*)=\log|L(s,\chi^*)|&\le C\mathcal L&&(\mathrm{Re}\,s>\tfrac78,\ |\mathrm{Im}\,s-t|\le\tfrac98),
\end{aligned}
\end{equation}
the second by convexity; only this upper bound for the real part is needed. Let $\mathcal C_r$ be the disc of radius $r$ about $2+it$; for $r<\frac98$ it lies in $\mathrm{Re}\,s>\frac78$. By \eqref{eq:logL}, the Borel--Carath\'{e}odory theorem from $\mathcal C_{9/8-\delta/4}$ to $\mathcal C_{9/8-\delta/2}$, and then the three-circles theorem with radii $\frac12<\frac98-\delta<\frac98-\frac\delta2$, as in \cite[Section 14.2]{Ti},
\[
\max_{\mathcal C_{9/8-\delta/2}}|\log L(s,\chi^*)|\ll_\delta\mathcal L,\qquad\max_{\mathcal C_{9/8-\delta}}|\log L(s,\chi^*)|\ll_\delta\mathcal L^{1-c},\qquad c=c(\delta)>0.
\]
The disc $\mathcal C_{9/8-\delta}$ contains $\sigma+it$ whenever $\frac78+\delta\le\sigma\le2$. Moreover $L(s,\chi)=L(s,\chi^*)\prod_{p\mid k}(1-\chi^*(p)p^{-s})$, and for $\mathrm{Re}\,s>0$ the absolute value of the product is at most $2^{\omega(k)}\ll_\eps k^\eps$. Hence
\begin{equation}\label{eq:lindelof-strip}
L(s,\chi)\ll_{\delta,\eps}(k(|t|+2))^\eps\qquad(\tfrac78+\delta\le\mathrm{Re}\,s\le2,\ t=\mathrm{Im}\,s).
\end{equation}
For $\mathrm{Re}\,s>0$ define $\check U(s)=\int_0^\infty U(x)x^{s-1}\,dx$. Integration by parts gives
\[
\check U(s)=-\frac1s\int_1^2U'(x)x^s\,dx,
\]
which provides the meromorphic continuation of $\check U$ to $\mathbb C$, with a single simple pole at $s=0$ of residue $1$, and $\check U(s)\ll_A|s|^{-A}$ for $|s|\ge1$ in vertical strips. In
\[
\sum_{m\ge1}\frac{\chi(m)}mU\Big(\frac mY\Big)=\frac1{2\pi i}\int_{(1)}L(1+s,\chi)\check U(s)Y^s\,ds
\]
move the line of integration to $\mathrm{Re}\,s=-\frac18+\delta$. Since $L(1+s,\chi)$ is entire, the only pole crossed is $s=0$, with residue $L(1,\chi)$, and by \eqref{eq:lindelof-strip} the new integral is $O_{\delta,\eps}(k^\eps Y^{-1/8+\delta})$. Take $\delta=\min\{\eps,\frac1{16}\}$.
\end{proof}

Now we can prove Theorem~\ref{thm:cond}.
\begin{proof}[Proof of Theorem~\ref{thm:cond}]
We make the changes of \cite[Section 6]{BK}; the hypothesis enters only through Lemma~\ref{lem:Lapprox}. Let $\alpha=\frac5{13}$, $H=R^{-3/5}$, $Q=R/d$ and $Y_d=(16Q^2)^\alpha$. In the proof of Proposition~\ref{prop:boundary}, truncate at $Y_d$ instead of $Z_d$, using Lemma~\ref{lem:Lapprox} for the characters $\chi_{\sigma a}$, whose moduli $k$ are at most $Z_d$; since $kY_d\le(16R^2)^{18/13}$, the factor $(kY_d)^\eps$ is $O(R^{3\eps})$. As in \cite[(6.3)--(6.7)]{BK}, with the second bound of Lemma~\ref{lem:combine} in place of \cite[(6.6)]{BK}, the error in Proposition~\ref{prop:boundary} becomes
\begin{align*}
&\ll R^\eps\big(R+R^{2-\alpha/4}H+R^{2-\alpha}H+R^{5/4+\alpha/2}H^{1/4}+R^{16/15+3\alpha/5}\big)\\
&=R^\eps\big(R+R^{339/260}+R^{66/65}+R^{84/65}+R^{253/195}\big).
\end{align*}
Apply Proposition~\ref{prop:expsum} with the exponent pair $(\kappa^*+\eta,\lambda^*+\eta)$ of \cite[(6.9)]{BK}, where, by \cite[Lemma 1.1]{TY},
\[
(\kappa^*,\lambda^*)=\Big(\frac{94971}{435968},\frac{530939}{871936}\Big),\qquad2\lambda^*-\kappa^*=1.
\]
With $p=\frac{1+\kappa^*}4=\frac{530939}{1743872}$, Proposition~\ref{prop:expsum} gives, since $2\lambda^*-\kappa^*=1$,
\[
\frac{V_M}M\ll(\rho M)^\eps\big(M^{1/4}+\rho^{p+\eta/4}M^{\eta/8}\big),
\]
and $\rho^{\eta/4}M^{\eta/8}\ll R^{\eta/2}$ for $\rho\asymp R$, $M\le R^2$. Taking $\eta$ small, \eqref{eq:Vhyp} holds with $(p_1,p_2)=(p,0)$, as in \cite[(6.11)]{BK}, and Lemma~\ref{lem:dual} bounds the dual sums by
\[
\ll R^\eps\big(RH^{-1/2}+R^{1+p}+R^{-1}H^{-3}\big)=R^\eps\big(R^{13/10}+R^{1+p}+R^{4/5}\big).
\]
Since $1+p=\frac{2274811}{1743872}>\frac{339}{260}$, Proposition~\ref{prop:Nsums} applies with $\vartheta=\frac{1+p}2$, and Corollary~\ref{cor:transfer} gives the theorem.
\end{proof}

\subsection{The Lindel\"{o}f hypothesis}\label{ssec:lindelof}

Now we prove Remark~\ref{rem:grh}, which says that if we assume the uniform generalized Lindel\"{o}f bound
\[
L(\tfrac12+it,\chi)\ll_\eps\big(q(1+|t|)\big)^\eps
\]
for primitive non-principal quadratic Dirichlet characters $\chi$ of conductor $q$, uniformly in $\chi$ and $t\in\mathbb R$, then the four formulas of Theorems~\ref{thm:neg} and~\ref{thm:pos} hold  with error terms changing to $O_\eps(N^{118/181+\eps})$.

\begin{proof}[Proof of the assertion in Remark~\ref{rem:grh}]
Assume that $L(\frac12+it,\chi)\ll_\eps(q(1+|t|))^\eps$ for every primitive non-principal Dirichlet quadratic character $\chi$ of conductor $q$ and every real $t$, uniformly in $\chi$ and $t$ (this follows from the Generalized Riemann Hypothesis. If $\chi$ modulo $k$ is induced by $\chi^*$, then $L(s,\chi)=L(s,\chi^*)\prod_{p\mid k}(1-\chi^*(p)p^{-s})$, and the absolute value of the product is at most $2^{\omega(k)}\ll_\eps k^\eps$ on $\mathrm{Re}\,s=\frac12$; so the bound holds for imprimitive $\chi$, with $k$ in place of $q$. The hypothesis enters only through the boundary terms, where it gives a square-root bound for the error term in the short sums of $L(1,\chi)$; this arithmetic observation is due to Heath-Brown \cite[Section 2]{HB}. For the dual sums we use the exponential-sum estimates of Section~\ref{sec:expsums}, with Bourgain’s exponent pair. Let $H=R^{-71/181}$ and $R\ge64$, so that $H\le\frac14$.

Under the hypotheses of Lemma~\ref{lem:PY},
\begin{equation}\label{eq:S-grh}
S^\sigma_\nu(Q^2;L,J)\ll_\eps R^\eps\min\big\{J^{1/2},\,LJ^{-1/2}+1\big\}\le R^\eps\big(L^{1/2}+1\big).
\end{equation}
For the first bound, let $m$ be odd, $m\ne\sq$, and $\psi_m(r)=(r/m)$, a non-principal character modulo $m$. For $n=\nu+8r$ and $8r_0\equiv\nu\pmod m$,
\[
\Big(\frac{\sigma(\nu+8r)}m\Big)=\Big(\frac{8\sigma}m\Big)\psi_m(r+r_0),
\]
so by the P\'{o}lya--Vinogradov inequality and partial summation the $n$-sum is $O(m^{1/2}\log2m)$, and its average over $m\asymp J$ is $O(J^{1/2}\log J)$. For the second, the proof of Lemma~\ref{lem:PY} with the assumed bound in place of the Weyl bound gives $O_\eps(n^\eps J^{1/2})$ for the $m$-sum when $\sigma n\ne\sq$, and the square value of $n$ (at most one, for $\sigma=1$) contributes $O(1)$.

By \eqref{eq:S-grh} and Remark~\ref{rem:general-char} with $(\beta,\gamma)=(0,\frac12)$, the boundary sums of Corollary~\ref{cor:boundary} equal their main terms up to
\[
\ll_\eps R^\eps\big(R+R^{3/2}H^{1/2}\big).
\]
By Remark~\ref{rem:pairs} with Bourgain's exponent pair, without the interpolation of Corollary~\ref{cor:Vbound}, \eqref{eq:Vhyp} holds with $(p_1,p_2)=(\frac{97}{336},\frac{13}{672})$, and Lemma~\ref{lem:dual} bounds the dual sums by
\[
\ll_\eps R^\eps\big(RH^{-1/2}+R^{433/336}H^{-13/336}+R^{-1}H^{-3}\big).
\]
For $H=R^{-h}$ the equation $\frac32-\frac h2=\frac{433}{336}+\frac{13h}{336}$ gives $h=\frac{71}{181}$, and then the exponents of $R$ in these error terms are
\[
1,\qquad\frac{236}{181},\qquad\frac{433}{362},\qquad\frac{236}{181},\qquad\frac{32}{181}.
\]
So Proposition~\ref{prop:Nsums} applies with $\vartheta=\frac{118}{181}$, and Corollary~\ref{cor:transfer} gives the assertion.
\end{proof}

 \subsection*{Acknowledgements}
 The authors would like to thank Prof. Biao Wang for his thoughtful and valuable suggestions during the preparation of this paper. The third author is very grateful to Prof. Genkai Zhang for his kind hospitality and encouragement.
 \subsection*{AI disclosure} The AI tools ChatGPT and Claude were used for discussing possible proof strategies, exploratory computations, editing and reorganization. All mathematical statements and proofs in the final manuscript
were written, reviewed, and verified by the authors, who take full responsibility for the
content of this manuscript.

\end{document}